\documentclass[12pt]{article}
\usepackage[utf8]{inputenc}
\usepackage[english]{babel}
\usepackage{amsmath,amssymb,amsthm}
\renewcommand\le{\leqslant}
\renewcommand\ge{\geqslant}

\newcommand\eps{\varepsilon}
\newcommand\R{\mathbb R}
\newtheorem*{theorem}{Theorem}
\newtheorem{lemma}{Lemma}
\newtheorem{statement}{Proposition}

\DeclareMathOperator{\codim}{codim}

\title{Product of octahedra is badly approximated in the $\ell_{2,1}$--metric}
\author{Yu.V. Malykhin, K.S. Ryutin}

\begin{document}

\maketitle

\begin{abstract}
    We prove that the cartesian product of octahedra 
    $B_{1,\infty}^{n,m}=B_1^n\times\ldots\times B_1^n$ ($m$ octahedra) is badly
    approximated by half--dimensional subspaces in mixed--norm:
    $d_{N/2}(B_{1,\infty}^{n,m},\ell_{2,1}^{n,m})\ge cm$, $N=mn$. As a corollary
    the orders for linear widths of H\"older--Nikolskii classes $H^r_p(\mathbb
    T^d)$ in the $L_q$ metric are obtained for $(p,q)$ in a certain set (a domain in the
    parameter space).
\end{abstract}

We consider the space $\R^N$, with $N=mn$.The set of coordinates $\{1,\ldots,N\}$ is split into   $m$ blocks $\Delta_1,\ldots,\Delta_m$ of cardinality $n$;
where $\Delta_s=\{(s-1)n+1,\ldots,sn\}$. For a vector  $x\in\R^N$
we denote by $x(i)$ its  $i$-th coordinate, and $x[s]$~--- its restriction on the 
$s$-block: $x[s]=(x(i))_{i\in\Delta_s}$.

We equip $\R^k$ with usual norms  $\|x\|_p=(\sum_{j=1}^k |x_j|^p)^{1/p},$ for
$1\le p <\infty$; $\|x\|_\infty=\max_{j=1}^k |x_j|$, and  $\R^{N}$~--- with mixed norms $$
\|x\|_{p,q} := \|x\|_{\ell_{p,q}^{n,m}} := \|y\|_{\ell_q^m},\quad\mbox{where
}y=(\|x[s]\|_{\ell_p^n})_{s=1}^m.
$$
Let $B_{p,q}^{n,m}$ be the unit ball in 
$\ell_{p,q}^{n,m}$ space, and  $B_p^k$, as usual, the unit ball of $\ell_p^k$.
So, $B_{1,\infty}^{n,m}$ is the cartesian product of  $m$ octahedra $B_1^n$. 
Throughout the paper we write $\|x\|_p:=\|x\|_{\ell_p^N}$ and $|x|:=\|x\|_2$ for brevity. 

Let us recall the definition of the Kolmogorov width of the subset  $M$ in the normed space $X$:
$$
d_n(M,X) = \inf_{\substack{L\subset X,\\\dim L\le n}} \sup_{x\in M}\inf_{y\in
L} \|x-y\|_X,
$$
and also the Gelfand width 
$$
d^n(M,X) = \inf_{\substack{L\subset X,\\\codim L\le n}} \sup_{x\in M\cap
L}\|x\|_X.
$$
We make use of the standard duality~\cite{IT,LGM}
$$
d_n(B_X,Y) = d^n(B_{Y^*},X^*),
$$
for  $X=(\R^N,\|\cdot\|_X)$ and $Y=(\R^N,\|\cdot\|_Y)$~---  finite dimensional
normed spaces with balls $B_X$ and $B_Y$, accordingly.

Throughout the paper with $c,c_1,c_2,\ldots,C,C_1,C_2,\ldots$ we denote different absolute positive constants (whose value may depend on the formula). The dependence on some parameters will be stated explicitly.

It is well known, that the calculation of widths of Sobolev classes is often
reduced to analogous calculation for  $B_p^N$ sets in $\ell_q^N$. While studying
widths of H\"older-Nikolskii classes for functions of several variables
E.M.~Galeev gave a lower estimate via widths of $B_{p,q}^{n,m}$.  In particular,
in~\cite{Gal90} he proved the following inequality:
$$
d_{N/2}(B_{1,\infty}^{n,m},\ell_{2,q}^{n,m})\ge c_qm^{1/q},\quad 1<q\le 2.
$$
(We certainly have simple upper estimate $d_k(B_{1,\infty}^{n,m},\ell_{2,q}^{n,m})\le
m^{1/q}$ for all $k$ and $1\le q\le \infty$.)
Galeev asked about the  $q=1$--case; A.D.~Izaak~\cite{Iz94} obtained the estimate 
\begin{equation}\label{izaak}
d_{N/2}(B_{1,\infty}^{n,m},\ell_{2,q}^{n,m})\ge cm\frac{\sqrt{\log\log m}}{\log
m}.
\end{equation}
Both proofs by Galeev and Izaak are based on the E.D. ~Gluskin~\cite{Gl87} method, that can not give the true order for $q=1$.

In this paper we establish the following result:
\begin{theorem}
    For all $m,n\in\mathbb N$ the following inequality holds
    $$
    d_{\lfloor mn/{2}\rfloor}(B_{1,\infty}^{n,m},\ell_{2,1}^{n,m})\ge cm,
    $$
    with $c>0$~--- an absolute constant.
\end{theorem}

There are some corollaries of the theorem. For example, for the width of the generalized octahedron  $V_m^N=B_\infty^N\cap mB_1^N$ we have an estimate
$d_{N/2}(V_m^N,\ell_q^N)\ge cm^{1/q}$, $q\in[1,2]$, with some absolute  $c$
(in Gluskin~\cite{Gl87} it depends on $q$), and also $$
cm^{1/q_2-1/p_2}\le d_{N/2}(B_{p_1,p_2}^{n,m},\ell_{q_1,q_2}^{n,m}) \le
m^{1/q_2-1/p_2},\quad p_1\le q_1\le 2,\; p_2\ge q_2.
$$

We obtain from theorem the following orders for the linear widths of the H\"older--Nikolskii classes $H^r_p=H^r_p(\mathbb T^d)$,
$0<r_1=\ldots=r_{l+1}<r_{l+2}=\ldots=r_d$ in  $L_q$ space, $1<p,q<\infty$
(for necessary definitions see e.g.~\cite{Gal90,Gal96}):
$$
\lambda_N(H^r_p,L_q)\asymp \begin{cases}
    \left(\frac{\log^l N}{N}\right)^{r_1-1/2+1/q}\log^{l/q}N,&\quad
    \frac1p+\frac1q<1,\;p\le 2,\;r_1>1-\frac1q,\\
    \left(\frac{\log^l N}{N}\right)^{r_1-1/p+1/q}\log^{l/q}N,&\quad
    2\le p<q,\;r_1>\frac1p-\frac1q.
\end{cases}
$$
Upper and lower bounds for  $\lambda_N(H^r_p,L_q)$ were obtained
in~\cite{Gal96,Iz96}, using the inequality~(\ref{izaak}) and, so, differed by a
power of $\log\log N$.  Substituting the estimate given by our Theorem into
Galeev's~\cite{Gal96} proof, we get the true order for   $\lambda_N(H^r_p,L_q)$.

It seems interesting to obtain from the  Theorem appropriate corollaries for the widths of Besov classes.

Let us recall that $N=mn$; throughout the paper we suppose that  $n$ and $m$ are large and $N$ is a multiple of $4$. In our proof we will use the well--known identity $d_k(B_1^N,\ell_2^N)=\sqrt{1-k/N}$. In dual terms for  $k=3N/4$ is means that
\begin{equation}\label{gen}
    \forall L\subset\R^N,\;\dim L\ge N/4\quad\Rightarrow\quad \exists x\in
    L\;\;|x|=1,\|x\|_\infty\ge 1/2.
\end{equation}

We will also need the next
\begin{lemma}
    Let $\mathcal X_1,\ldots,\mathcal X_m$ be finite sets of vectors in 
    $\R^d$, with $\sum_{x\in\mathcal X_s}|x|^2\le 1$ for all $s$. Then there exists some vector  $\eps\in\{0,1,-1\}^d$, $\|\eps\|_1\ge cd$, such that
    $$
    \max_{s=1,\ldots,m}\left(\sum_{x\in\mathcal X_s}\langle
    \eps,x\rangle^2\right)^{1/2} \le C\log^{1/2}\left(\frac{m}{d}+2\right).
    $$
\end{lemma}
This lemma generalizes Theorem 3 by Gluskin~\cite{Gl88}, which dealt with vectors  $x_1,\ldots,x_m$ (i.e. single element sets  $\mathcal X_s = \{x_s\}$).
Gluskin strengthens and develops  the preceding result by B.S. Kashin~\cite{K85};
the approaches to vector balancing problems introduced in these papers were later improved by several authors.  The proof of lemma  1 from known results on vector balancing and properties of gaussian measure on $\R^d$ will be given at the end of the paper. 


\begin{proof}[The proof of the theorem]
By duality,
$d_{N/2}(B_{1,\infty}^{n,m},\ell_{2,1}^{n,m})=d^{N/2}(B_{2,\infty}^{n,m},\ell_{\infty,1}^{n,m})$.
We have to find in any subspace  $L\subset\R^N$ of the half--dimension ($N/2$) a vector $x$, with  $\|x\|_{\infty,1}\ge cm$,
$\|x\|_{2,\infty}\le C$. These requirements will be fulfilled provided
\begin{itemize}
    \item $\ell_2$-norm of any block is bounded: $|x[s]|\le C$;
    \item there are  $\ge cm$ blocks $s$, such that $|x(i_s)|\ge 1/4$ for some  $i_s\in\Delta_s$.
\end{itemize}

The required vector will be obtained as a certain linear combination of vectors 
$x_1,\ldots,x_l\in L$, where subspaces $Z_j\subset Z_{j-1}\subset L$ and
vector $x_j\in Z_j$ will be produced on the  $j$-th step of the  construction desribed below.
On the first step  we let $Z_1:=L$ and take any vector $x_1$ given by ~(\ref{gen}). We have $x_1\in L$, $|x_1|=1$ and $\|x_1\|_\infty\ge 1/2$.
So, $|x_1(i_1^*)|\ge 1/2$ for some coordinate $i_1^*$; we call this coordinate  ``large''; if there are several such coordinates  (certainly, no more than $4$), we choose  only one of them. We denote by $i_j^*$ the large coordinate produced on the  $j$-th step.

The numbers $(x_k(i))_{\substack{1\le k\le j\\1\le
i\le N}}$ obtained after the first $j$ steps  are written into an array where $k$ is the row number  and  $i$~--- the column number. 

After $j$-th step we introduce the set of  ``vanishing'' coordinates 
$\Lambda_j\subset\{1,\ldots,N\}$. This set consist of:
\begin{itemize}
    \item[(i)] coordinates  $i$, with a large sum of squares (in fixed coordinate):
        $x_1(i)^2+\ldots+x_j(i)^2\ge 1/n$;
    \item[(ii)] all coordinates of all blocks $s$, that contain some ``large''
        coordinate: $i_k^*\in\Delta_s$, $1\le k\le j$;
    \item[(iii)] all coordinates from any block $s$, with a large sum of squares (in block):
        $|x_1[s]|^2 + \ldots + |x_j[s]|^2 \ge 1$.
\end{itemize}
At step $(j+1)$ we apply (\ref{gen}) to the subspace  $Z_{j+1}:=\{x\in L\colon
x(i)=0\;\forall i\in\Lambda_j\}$ and obtain a vector $x_{j+1}$.

We can apply~(\ref{gen}) if
$|\Lambda_j|\le N/4$.  Let us show that this holds for $j\le cm$.  Really, at
each step only one new block with a large coordinate is introduced, therefore we
get no more than   $jn\le cN<N/12$ vanishing coordinates. Then, since the sum
(on all $k\le j$ and   $i\le N$) of squares of numbers $x_k(i)$ does not exceed
$j$, we have not more than $j/(1/n)=jn<N/12$ columns with sum of squares
$\ge 1/n$. Similarly we can upper estimate the number of vanishing coordinates
from blocks with large sum of squares (in block).

So, we will make $l$ steps, $l\asymp m$ (i.e. $c_1m\le l\le c_2m$) and will obtain vectors $x_1,\ldots,x_l$.
All numbers $x_j(i)$ can be divided into three groups:
\begin{itemize}
    \item ``large'': $x_j(i_j^*)$;
    \item ``intermediate'':  $x_j(i)$, such that $\sum_{k\le j}x_k(i)^2\ge 1/n$;
    \item ``small'': all other numbers.
\end{itemize}
We note that (according to our construction) if some number $x_j(i)$ is large or intermediate then on all steps after $j$-th  the $i$-th coordinate vanishes: $x_{j+1}(i)=\ldots=x_l(i)=0$.
Let  $x_j$ be written as  $v_j+w_j$, where $v_j$ consist of all large and intermediate numbers  and  $w_j$ of small.

We will obtain the required vector  $x$ as a (balanced) sum 
$$
x=\sum_{j=1}^l \eps_j x_j = \sum_{j=1}^l \eps_j v_j + \sum_{j=1}^l \eps_j w_j =
v + w,
$$
where  $\eps=(\eps_j)_{j=1}^l\in\{0,1,-1\}^l$ will be given by lemma. It is clear that $x\in L$. We have
$|x[s]|\le |v[s]|+|w[s]|$. 
Since for any  $i$ the column  $(v_j(i))_{j=1}^l$ contains only one nonzero number we get
$$
|v[s]|^2 = \sum_{i\in\Delta_s} v(i)^2 = \sum_{i\in\Delta_s} (\sum_{j=1}^l
\eps_j v_j(i))^2 \le \sum_{i\in\Delta_s}\sum_{j=1}^l v_j(i)^2 =
\sum_{j=1}^l|v_j[s]|^2 \le C.
$$
The last inequality holds because our construction and  (iii) imply even an estimate $\sum_{j=1}^l|x_j[s]|^2 \le C$.

It is clear that  $|v(i_j^*)|=|v_j(i_j^*)|\ge 1/2$,
if $|\eps_j|=1$. All large coordinates are from different blocks by our construction, so 
$\|v\|_{\infty,1}\ge \frac12\|\eps\|_1$.
In order to finish the proof we need to find vector  $\eps$, $\|\eps\|_1\ge cm$, such that 1)
$\ell_2$-norm of $w$ in any block is bounded by some absolute constant 2) $|w(i)|<1/4$ for all $i\in
I:=\{i_1^*,\ldots,i_l^*\}$.

Let us consider columns 
$f_i := (w_j(i))_{j=1}^l\in\R^l$; we have $|f_i|\le1/\sqrt{n}$ because of (i).
We apply the lemma to the following sets of vectors: 
for each block  $s=1,\ldots,m$ we introduce the set $\{f_i\colon i\in\Delta_s\}$, and we also add single--point  sets  $\{f_i/|f_i|\}$ for $i\in I$.
Lemma gives us some vector  $\eps$, such that
$|w[s]|^2 = \sum_{i\in\Delta_s}\langle \eps,f_i\rangle^2 \le C$,
and also $|w(i)| = |\langle\eps,f_i\rangle| \le C|f_i|\le C/\sqrt{n} < 1/4$ for
$i\in I$. The proof is finished.

\end{proof}

Let us prove lemma  1.
Let  $\xi=(\xi_1,\ldots,\xi_d)$ denote the standard gaussian vector in  $\R^d$,
$\gamma_d$~--- be the gaussian measure in $\R^d$: $\gamma_d(A)=\mathsf P(\xi\in
A)=\int_A(2\pi)^{-d/2}\exp(-|x|^2/2)\,dx$, $A\subset\R^d$. Let
$\Psi(t) := \gamma_1(-t,t)$.
We recall well-known estimates~\cite{Gl88},\cite{lifsh}: $\Psi(t)\ge 1-e^{-t^2/2}\ge
\exp(-2\exp(-t^2/2))$, $t\ge 1$.

We will use the following results. 

\begin{theorem}[$S$-inequality]
Let $K\subset \R^d$ be a convex closed centrally--symmetric set, $w\in\R^d$, and $P=\{x:|\langle w,x\rangle|\le 1\}$ --- the strip of the same gaussian measure $\gamma_n(K)=\gamma_n(P)$. Then the inequality 
$\gamma_n(tK)\ge \gamma_n(tP)$ holds for $t\ge 1$ and the opposite inequality
holds for $t\in(0,1)$.
\end{theorem}
The proof was obtained by  R.~Latala and K.~Oleszkiewicz~\cite{latol}.

\begin{theorem}[Gaussian Correlation Conjecture]
For any convex centrally--symmetric sets $K_1,K_2\subset \R^d$ holds 
$\gamma_d(K_1\cap K_2)\ge \gamma_d(K_1)\gamma_d(K_2)$.
\end{theorem}
This theorem was proved in a recent paper by  T.~Royen~\cite{roy}; its
particular case when $K_1$ is an ellipsoid (it was proved in \cite{H99}) is
sufficient for our purposes.

Besides it, we will need the following particular case of lemma 3.2 from~\cite{Gi97}:
\begin{statement}
    For any convex centrally--symmetric set  $V\subset\R^d$,
    such that $\gamma_d(V)\ge 2^{-d/7}$, and any vectors  $u_1,\ldots,u_d\in B_2^d$ there exists vector  
    $\eps\in\{0,1,-1\}^d$, $\|\eps\|_1\ge d/2$, such that    $\sum_{i=1}^d\eps_iu_i\in 4V$.
\end{statement}

We will use it when  $\{u_1,\ldots,u_d\}$ are standard basis vectors in  $\R^d$, and get a vector  $\eps\in 4V$.  Let
$$
h_s(v):=(\sum_{x\in\mathcal X_s}\langle v,x\rangle^2)^{1/2},\quad
\mathcal E_s:=\{v\in\R^d\colon h_s(v)\le 1\}.
$$
For a standard gaussian vector $\mathsf{E}h_s(\xi)^2\le
\sum_{x\in\mathcal X_s}|x|^2 \le 1$, so we get the estimate
\begin{equation}\label{ebound}
\gamma_d(2\mathcal E_s) = 1-\mathsf P(h_s(\xi)>2)\ge
1-\frac{\mathsf{E}h_s(\xi)^2}{4}\ge 3/4 > \Psi(1).
\end{equation}
By $S$--inequality, using~(\ref{ebound}), we get 
$$
\gamma_d(2t\mathcal E_s) \ge \Psi(t) \ge 
\exp(-2\exp(-t^2/2)),\quad t\ge 1.
$$

Let $V=\bigcap_{s=1}^m(2t\mathcal E_s)$;
applying the  Gaussian Correlation Conjecture for ellipsoids, we obtain
$$
\gamma_d(V)\ge\prod_{s=1}^m\gamma_d(2t\mathcal E_s) \ge 
\exp(-2m\exp(-t^2/2)).
$$
For $t\asymp \log^{1/2}(m/d+2)$ we have $\gamma_d(V)\ge 2^{-d/7}$.  Therefore,
we can apply proposition 1 to $V$, which gives us the required vector $\eps$.

\end{document}